\documentclass{amsart}

\usepackage{latexsym,enumerate}
\usepackage{amsmath,amsthm,amsopn,amstext,amscd,amsfonts,amssymb}
\usepackage[ansinew]{inputenc}
\usepackage{verbatim}
\usepackage{graphicx}
\usepackage{pstricks}
\usepackage{wasysym}
\usepackage[all]{xy}
\usepackage{epsfig} 
\usepackage{graphicx,psfrag} 
\usepackage{subfigure}
\usepackage{pstricks,pst-node}

\usepackage{multicol}
\usepackage{color}
\usepackage{colortbl}

\newcommand{\RR}{\mathbb{R}}
\newcommand{\ZZ}{\mathbb{Z}}

\newtheorem{theorem}{Theorem}[section]
\newtheorem{lemma}[theorem]{Lemma}
\newtheorem{proposition}[theorem]{Proposition}

\newtheorem{definition}[theorem]{Definition}

\newtheorem{remark}[theorem]{Remark}

\DeclareMathOperator{\Vol}{Vol}

\newcommand{\spb}[1]{\smallskip}
\newcommand{\mpb}[1]{\medskip}
\newcommand{\bpb}[1]{\bigskip}

\newcommand{\p}{\partial}

\renewcommand{\a}{\alpha}

\renewcommand{\d}{\delta}

\newcommand{\G}{\Gamma}

\renewcommand{\l}{\lambda}

\newcommand{\s}{\sigma}

\begin{document}
\DeclareGraphicsExtensions{.jpg,.pdf,.mps,.png}

\title{A note on isoperimetric inequalities of Gromov hyperbolic manifolds and graphs}

\author[\'{A}lvaro Mart\'{\i}nez-P\'erez]{\'{A}lvaro Mart\'{\i}nez-P\'erez}

\author[Jos\'e M. Rodr{\'\i}guez]{Jos\'e M. Rodr{\'\i}guez}

\date{\today}

\begin{abstract}
In this paper we study the relationship of hyperbolicity and (Cheeger) isoperimetric inequality in the context of Riemannian manifolds and graphs.
We characterize the hyperbolic manifolds and graphs (with bounded local geometry) verifying this isoperimetric inequality, in terms of their Gromov boundary improving similar results from a previous work.
In particular, we prove that having a pole is a necessary condition and, therefore, it can be removed as hypothesis.
\end{abstract}

\maketitle{}

\

\


\


\

\

Author names and affiliations:

\bigskip

\'{A}lvaro Mart\'{\i}nez-P\'erez,
Facultad CC. Sociales de Talavera,
Avda. Real Fábrica de Seda, s/n. 45600 Talavera de la Reina, Toledo, Spain,

\noindent
alvaro.martinezperez$@$uclm.es

\noindent
ORCID 0000-0002-1344-6189

\bigskip

Jos\'e M. Rodr\'{\i}guez (corresponding author),
Departamento de Matem\'aticas, Universidad Carlos III de Madrid,
Avenida de la Universidad 30, 28911 Legan\'es, Madrid, Spain, Phone number 34 91 624 9098,

\noindent
jomaro$@$math.uc3m.es

\noindent
ORCID 0000-0003-2851-7442

\newpage

\section{Introduction}

Isoperimetric inequalities are of interest in pure and applied mathematics (see, e.g., \cite{C2}, \cite{Po}).
There are close connections between isoperimetric inequality and some conformal invariants of Riemannian manifolds and graphs, namely
Poincar\'e-Sobolev inequalities,
the bottom of the spectrum of the Laplace-Beltrami operator, the exponent of convergence,
and the Hausdorff dimensions of the sets of both bounded geodesics and escaping geodesics in a negatively curved surface
(see \cite{APR}, \cite{BJ}, \cite[p.228]{Bu}, \cite{Ch}, \cite{FM1}, \cite{FM2}, \cite{FMP1}, \cite{FMP2}, \cite{FR1}, \cite{MRT}, \cite{P}, \cite[p.333]{S}).
The Cheeger isoperimetric inequality is closely related to the project of Ancona
on the space of positive harmonic functions of Gromov-hyperbolic manifolds and graphs (\cite{A1}, \cite{A2} and \cite{A3}).
In fact, in the study of the Laplace operator on a hyperbolic manifold or graph $X$,
Ancona obtained in those three last papers interesting results, under the additional assumption that the bottom of the spectrum of the Laplace spectrum $\l_1(X)$ is positive.
The well-known Cheeger inequality $\l_1(X) \ge \frac14 \,h(X)^2$, where $h(X)$ is the isoperimetric constant of $X$,
guarantees that $\l_1(X)>0$ when $h(X) > 0$ (see \cite{Bu1} for a converse inequality).
Hence, the results of this paper are useful in order to apply these Ancona's results.

There is a natural connection between hyperbolicity and Cheeger isoperimetric inequality. In fact, one of the alternative definitions of Gromov hyperbolicity uses some kind of isoperimetric inequality (see \cite{ABCD}, \cite{G1}).

Cao proved in \cite{Cao} that hyperbolicity with an extra hypothesis on the Gromov boundary implies (Cheeger) isoperimetric inequality
(an extra hypothesis is necessary, since there exist hyperbolic graphs without isoperimetric inequality, as the Cayley graph of the group $\ZZ$).

In \cite{MR1} we studied the relationship of hyperbolicity and Cheeger isoperimetric inequality in the context of Riemannian manifolds and graphs
with bounded local geometry.

Given any Riemannian $n$-manifold $M$, the Cheeger isoperimetric constant of $M$ is defined as
\[h(M) = \inf_A \frac{\Vol_{n-1}(\partial A)}{\Vol_{n}(A)} ,\]
where $A$ ranges over all non-empty bounded open subsets of $M$,
and $\Vol_{k}(B)$ denotes the $k$-dimensional Riemannian volume of the set $B$.

Given any graph $\G=(V,E)=(V(\G),E(\G))$, let us consider the natural length metric $d_{\,\Gamma}$ where every edge has length 1.
For any graph $\Gamma$, any vertex $v\in V$ and any $k\in \mathbb{N}$, let $S(v,k):=\{w\in V \, | \, d_{\,\Gamma}(v,w)=k \}$.
As usual, we denote by $B(v,k)$ and $\bar{B}(v,k)$ the open and closed balls, respectively.

The combinatorial Cheeger isoperimetric constant of $\Gamma$ is defined to be
\[h(\Gamma) = \inf_A \frac{|\partial A|}{|A|} ,\]
where $A$ ranges over all non-empty finite subsets of vertices in $\Gamma$,
$\partial A = \{v \in \Gamma \, | \, d_{\,\Gamma}(v,A) = 1\}$ and $|A|$ denotes the cardinality of $A$.

A Riemannian manifold or graph $X$ satisfies the (Cheeger) \emph{isoperimetric inequality} if $h(X)>0$, since in this case
$$
\Vol_{n}(A) \le h(X)^{-1} \Vol_{n-1}(\partial A) ,
$$
for every bounded open set $A \subseteq X$ if $X$ is a Riemannian $n$-manifold, and
$$
|A| \le h(X)^{-1} |\partial A| ,
$$
for every finite set $A \subseteq V(X)$ if $X$ is a graph.

Along the paper, we just consider manifolds and graphs $X$ which are connected.
This is not a loss of generality, since if $X$ has connected components $\{X_j\}$, then $h(X)=\inf_j h(X_j)$.

Let $(X,d_X)$ and $(Y,d_Y)$  be two metric spaces. A map $f: X\longrightarrow Y$ is said to be
an $(\alpha, \beta)$-\emph{quasi-isometric embedding}, with constants $\alpha\geq 1,\
\beta\geq 0$, if for every $x, y\in X$:
$$
\alpha^{-1}d_X(x,y)-\beta\leq d_{\,Y}(f(x),f(y))\leq \alpha \, d_X(x,y)+\beta.
$$
The function $f$ is $\varepsilon$-\emph{full} if
for each $y \in Y$ there exists $x\in X$ with $d_Y(f(x),y)\leq \varepsilon$.

A map $f: X\longrightarrow Y$ is said to be
a \emph{quasi-isometry}, if there exist constants $\alpha\geq 1,\
\beta,\varepsilon \geq 0$ such that $f$ is an $\varepsilon$-full
$(\alpha, \beta)$-quasi-isometric embedding.
Two metric spaces $X$ and $Y$ are \emph{quasi-isometric} if there exists
a quasi-isometry $f:X\longrightarrow Y$.
One can check that to be quasi-isometric is an equivalence relation.

A graph $\Gamma$ is said to be $\mu$-\emph{uniform} if each vertex $p$ of $V$ has at most $\mu$ neighbors, i.e.,
\[\sup\big\{|N(p)| \,  \big| \,\, p\in V(\G)\big\}\leq \mu.  \]
If a graph $\Gamma$ is $\mu$-uniform for some constant $\mu$ we say that $\Gamma$ is \emph{uniform}
or that it has \emph{bounded local geometry}.

A Riemannian $n$-manifold $M$ has \emph{bounded local geometry} if there exist positive constants $r,c,$ such that for every $x \in M$ there is a diffeomorphism $F : B(x,r) \rightarrow \RR^n$
with
$$
\frac{1}{c}\,d(x_1,x_2) \le \| F(x_1) - F(x_2) \| \le c\,d(x_1,x_2)
$$
for every $x_1,x_2 \in B(x,r)$.

The \emph{injectivity radius} inj$(x)$ of \emph{$x\in M$} is defined as the supremum of those $r>0$ such that $B(x,r)$ is simply connected or, equivalently,
as half the infimum of the lengths of the (homotopically non-trivial) loops based at $x$.
The \emph{injectivity radius} inj$(M)$ \emph{of $M$} is the infimum over $x\in M$ of inj$(x)$.

\begin{remark} \label{r:Ricci}
If $M$ has positive injectivity radius and a lower bound on its Ricci curvature, then $M$ has bounded local geometry \cite{AC}.
\end{remark}

A celebrated theorem of Kanai in \cite{K} states that quasi-isometries preserve isoperimetric inequalities
between Riemannian manifolds and graphs with bounded local geometry.
This result also holds with weaker hypotheses in the context of Riemann surfaces \cite{CGPR}, \cite{GPPRT}.

\medskip

Let $X$ be a metric space. Fix a base point $o\in X$ and for
$x,x'\in X$ let
$$(x|x')_o=\frac{1}{2}\big(d(x,o)+d(x',o)-d(x,x')\big).$$
The number $(x|x')_o$ is non-negative and it is called the \emph{Gromov product} of $x,x'$ with respect to $o$.

\begin{definition} A metric space $X$ is \emph{(Gromov)
hyperbolic} if it satisfies the $\delta$-inequality
\[(x|y)_o\geq \min\{(x|z)_o,(z|y)_o\}-\delta\] for some $\delta\geq
0$, for every base point $o\in X$ and all $x,y,z \in X$.
\end{definition}

We denote by $\d(X)$ the sharp hyperbolicity constant of $X$:
$$
\d(X)= \sup \Big\{ \min\{(x|z)_o,(z|y)_o\} - (x|y)_o \,\big| \;\, x,y,z,o \in X \Big\}.
$$
Hence, $X$ is hyperbolic if and only if $\d(X)<\infty$.

The theory of Gromov hyperbolic spaces was introduced by M. Gromov for
the study of finitely generated groups (see \cite{G1}).
The concept of Gromov hyperbolicity grasps the essence of negatively curved
spaces like the classical hyperbolic space, Riemannian manifolds of
negative sectional curvature bounded away from $0$, and of discrete spaces like trees
and the Cayley graphs of many finitely generated groups. It is remarkable
that a simple concept leads to such a rich
general theory (see \cite{ABCD, GH, G1}).
This theory has been developed from a geometric point of view to the extent of making
hyperbolic spaces an important class of metric spaces to be studied on their
own (see, e.g., \cite{BH,BBI,BS,GH,V}). In the last years, Gromov hyperbolicity
has been intensely studied in graphs (see, e.g.,
\cite{BRS,
BHB1,K21,K22,M,M2,PRSV,RT2,Sha1,Sha2,Sha3,WZ} and the references therein).
Gromov hyperbolicity, specially in graphs, has found applications in different areas such as phylogenetics (see
\cite{DHHKMW,DMT}), real networks (see \cite{AAD,ASM,CMN,KPKVB,MoSoVi}) or the secure transmission
of information and virus propagation on networks (see \cite{K21,K22}).

We want to remark that the main examples of hyperbolic graphs are the trees.
In fact, the hyperbolicity constant of a metric space can be viewed as a measure of
how ``tree-like'' the space is, since those spaces $X$ with $\delta(X) = 0$ are precisely the metric trees.
This is an interesting subject since, in
many applications, one finds that the borderline between tractable and intractable
cases may be the tree-like degree of the structure to be dealt with
(see, e.g., \cite{CYY}).

\smallskip

In \cite{MR1} we characterized, in terms of their Gromov boundary, the uniform hyperbolic graphs  and a large class of hyperbolic manifolds satisfying isoperimetric inequality (see, respectively, Theorems \ref{t:iigraph} and \ref{th: isoperimetric_manifolds}).
In Theorems \ref{t:iigraph}, \ref{th: isoperimetric_manifolds} and \cite[Theorem 1.1]{Cao}
it is used the hypothesis of the existence of a pole
(in fact, although \cite[Theorem 1.1]{Cao} apparently uses a different hypothesis, actually, in hyperbolic spaces, it is equivalent
to the existence of a pole).
The hypothesis on Gromov hyperbolicity is natural (we need it in order to deal with the Gromov boundary).
However, although the hypothesis on the pole is technically needed in the proofs, it does not look natural.
The goal of this paper is to remove this hypothesis in the statements of Theorems \ref{t:iigraph} and \ref{th: isoperimetric_manifolds}. Thus, in Theorems \ref{t:graphpole} and
\ref{th:manifoldspole} we prove that having a pole is not needed as an hypothesis but it is also a necessary condition of the characterization.

\section{Some previous results}

Recall that a \emph{geodesic space} is a metric space such that for every couple of points there exists a geodesic joining them.

\begin{definition}  A geodesic space $X$ has a \emph{pole} in a point $v$ if there
exists $M > 0$ such that each point of $X$ lies in an $M$-neighborhood of some geodesic
ray emanating from $v$.
\end{definition}

If $X$ is a geodesic metric space and $x_1,x_2,x_3\in X$, the union
of three geodesics $[x_1 x_2]$, $[x_2 x_3]$ and $[x_3 x_1]$ is a
\emph{geodesic triangle} that will be denoted by $T=\{x_1,x_2,x_3\}$
and we will say that $x_1,x_2$ and $x_3$ are the vertices of $T$.
We say that $T$ is $\d$-{\it thin} if any side of $T$ is contained in the
$\d$-neighborhood of the union of the two other sides.
We denote by $\d_{th}(T)$ the sharp thin constant of $T$, i.e., $ \d_{th}(T):=\inf\{\d\ge 0:
\, T \, \text{ is $\d$-thin}\,\}. $ The space $X$ is $\d$-thin (or satisfies the {\it Rips condition} with
constant $\d$) if every geodesic triangle in $X$ is $\d$-thin.
We denote by $\d_{th}(X)$ the sharp thin constant of $X$, i.e., $ \d_{th}(X):=\sup\{\d_{th}(T):
\, T \, \text{ is a geodesic triangle in $X$}\,\}. $

It is well-known that a geodesic metric space is hyperbolic if and only if it satisfies the Rips condition for some constant (see, e.g., \cite{ABCD,GH}).
In the classical references on this subject (see, e.g., \cite{ABCD,GH})
appear many different definitions of Gromov hyperbolicity, which are equivalent in the sense
that if $X$ is $\d$-hyperbolic with respect to one definition,
then it is $\d'$-hyperbolic with respect to another definition (for some $\d'$ related to~$\d$).

In order to consider a graph $G$ as a geodesic metric space, identify (by an isometry)
any edge $uv\in E(G)$ with the interval $[0,1]$ in the real line;
then the edge $uv$ (considered as a graph with just one edge)
is isometric to the interval $[0,1]$.
Thus, the points in $G$ are the vertices and, also, the points in the interior
of any edge of $G$.
In this way, any connected graph $G$ has a natural distance
defined on its points, induced by taking shortest paths in $G$,
and we can see $G$ as a metric graph.

\smallskip

Let us adapt the following definition from \cite{BS} where we introduce the constant $\varepsilon_0$ for convenience. Notice that for bounded metric spaces both definitions coincide. Since herein this property will be always applied to compact spaces all the results work as well with the original definition.

\begin{definition} Given a metric space $(X,d)$ and a constant $S>1$, we say that $(X,d)$ is \emph{$S$-uniformly perfect} if there exists some $\varepsilon_0>0$ such that for every $x\in X$ and every $0<\varepsilon \le \varepsilon_0$ there exist a point $y\in X$ such that
$\frac{\varepsilon}{S} <d(x,y) \le \varepsilon$. We say that $(X,d)$ is \emph{uniformly perfect} if there exists some $S$ such that $(X,d)$ is $S$-uniformly perfect.
\end{definition}

\section{Hyperbolic graphs}

Let us recall the concepts of geodesic and sequential boundary of a hyperbolic space and some basic properties. For further information and proofs we refer the reader to \cite{BH,BS,GH,G1}.

\smallskip

Let $X$ be a hyperbolic space and $o\in X$ a base point.

The \emph{relative geodesic boundary} of $X$ with respect to the base-point $o$ is
\[\partial_o^g X := \{ [\gamma]  \, | \, \gamma: [0,\infty) \to X \mbox{ is a geodesic ray with } \gamma(0) = o\},\]
where  $\gamma_1 \sim \gamma_2$ if there exists some $K>0$ such that
$d(\gamma_1(t),\gamma_2(t))<K$, for every $t\geq 0.$

In fact, the definition above is independent from the base point.
Therefore, the set of classes of geodesic rays is called \emph{geodesic boundary} of $X$, $\partial^gX$.  Herein, we do not distinguish between the geodesic ray and its image.

A sequence of points $\{x_i\}\subset X$ \emph{converges to infinity}
if \[\lim_{i,j\to \infty} (x_i|x_j)_o=\infty.\] This property is
independent of the choice of $o$ since
\[|(x|x')_o-(x|x')_{o'}|\leq d(o,o')\] for any $x,x',o,o' \in X$.

Two sequences $\{x_i\},\{x'_i\}$ that converge to infinity are
\emph{equivalent} if \[\lim_{i\to \infty} (x_i|x'_i)_o=\infty.\]
Using the $\delta$-inequality, we easily see that this defines an
equivalence relation for sequences in $X$ converging to infinity.
The \emph{sequential boundary at infinity} $\partial_\infty X$ of $X$ is
defined to be the set of equivalence classes of sequences
converging to infinity.

Note that given a geodesic ray $\gamma$, the sequence $\{\gamma(n)\}$ converges to infinity and two equivalent rays induce equivalent sequences.
Thus, in general, $\partial^g X\subseteq \partial_\infty X$.

We say that a metric space is \emph{proper} if every closed ball is compact.
Every uniform graph and every complete Riemannian manifold are proper geodesic spaces.

\begin{proposition}\cite[Chapter III.H, Proposition 3.1]{BH}\label{Prop: equiv_boundary} If $X$ is a proper hyperbolic geodesic space, then the natural map from $\partial^g X$ to $\partial_\infty X$ is a bijection.
\end{proposition}

For every $\xi, \xi' \in \partial_\infty X$, its Gromov product
with respect to the base point
$o\in X$ is defined as
\[ (\xi|\xi')_o =  \inf \ \liminf_{i\to \infty} (x_i|x'_i)_o,\]
where the infimum is taken over all sequences $\{x_i\} \in \xi $, $\{x'_i\} \in \xi' $.

A metric $d$ on the sequential boundary at infinity
$\partial_\infty X$ of $X$ is said to be \emph{visual}, if there are $o\in X$, $a > 1$ and positive
constants $c_1$, $c_2$, such that
\[c_1a^{-(\xi|\xi')_o} \leq  d(\xi, \xi') \leq c_2a^{-(\xi|\xi')_o}\]
for all $\xi, \xi' \in \partial_\infty X$. In this case, we say that $d$ is a visual metric
with respect to the base point $o$ and the parameter $a$.

\begin{theorem}\cite[Theorem 2.2.7]{BS} Let $X$ be a hyperbolic space. Then for any $o \in X$,
there is $a_0 > 1$ such that for every $a \in (1, a_0]$ there exists a metric $d$ on
$\partial_\infty X$, which is visual with respect to $o$ and $a$.
\end{theorem}

\begin{remark} Notice that for any visual metric, $\partial_\infty X$ is bounded and complete.
\end{remark}

\begin{theorem} \label{t:iigraph}
Given a hyperbolic uniform infinite graph $\Gamma$ with a pole, then $h(\Gamma)>0$ if and only if $\partial_\infty \Gamma$ is uniformly perfect for some visual metric.
\end{theorem}

\section{Hyperbolic manifolds}

Let us recall the following definitions from \cite{K}.

A subset $A$ in a metric space $(X,d)$ is called \emph{$r$-separated},
$r>0$, if $d(a,a')\geq r$ for any distinct $a,a'\in A$. Note that
if $A$ is maximal with this property, then the union $\cup_{a\in
A} B_r(a)$ covers $X$. A maximal $r$-separated set $A$ in a metric
space $X$ is called an $r$-\emph{approximation} of $X$.

Let $X$ be a complete Riemannian manifold and denote by $d$ the induced metric.
Given any $\varepsilon$-approximation $A_\varepsilon$ of $X$,
the graph $\Gamma_{A_\varepsilon}=(V,E)$ with $V=A_\varepsilon$ and $E:=\{xy \, | \, x,y \in A_\varepsilon \mbox{ with } 0<d(x,y)\leq 2\varepsilon\}$ is called an $\varepsilon$-\emph{net}.

\begin{proposition}\cite[Lemma 4.5]{K} \label{Prop: Kanai_net} Suppose that $X$ is a complete Riemannian manifold with bounded local geometry and let $\Gamma$ be an $\varepsilon$-net in $X$.
Then, $h(X)>0$ if and only if $h(\Gamma)>0$.
\end{proposition}

Note that the results in \cite{K} require $M$ to have positive injectivity radius and a lower bound on its Ricci curvature instead of bounded local geometry,
but the proofs in \cite{K} just use that there are uniform lower and upper bounds for the volume of the balls $B(x,r)$ which do not depend on $x\in M$ for $0<r<r_0$
(and we have these uniform bounds with bounded local geometry).
Hence, the results in \cite{K} also hold with the weaker hypothesis of bounded local geometry.

\begin{proposition}\cite[Lemma 2.5]{K} \label{Prop: net_qi}
Suppose that $X$ is a complete Riemannian manifold with bounded local geometry and let $\Gamma$ be an $\varepsilon$-net in $X$.
Then, $X$ and $\Gamma$ are quasi-isometric.
\end{proposition}

\begin{theorem}\cite[p.88]{GH}\label{th: stability_hyp}
If $f:X \rightarrow Y$ is a quasi-isometry between geodesic spaces, then $X$ is hyperbolic if and only if $Y$ is hyperbolic.
\end{theorem}

\begin{proposition}\cite[Proposition 5.6]{MR1}\label{Prop: qi-pole}
Suppose $X,Y$ are proper hyperbolic geodesic spaces and $f:X \to Y$ is a quasi-isometry.
If $X$ has a pole in $v$, then $Y$ has a pole in $f(v)$.
\end{proposition}

\begin{lemma}\cite[Lemma 2.3]{K} \label{lema: uniform-net} Every $\varepsilon$-net in a complete Riemannian manifold with bounded local geometry is uniform.
\end{lemma}

\begin{theorem}\cite[Theorem 5.12]{MR1}\label{th: isoperimetric_manifolds}
Let $X$ be a non-compact complete Riemannian manifold with bounded local geometry.
Assume that $X$ is hyperbolic and has a pole.
Then, $h(X)>0$ if and only if $\partial_\infty X$ is uniformly perfect.
\end{theorem}

In \cite{Cao} appear several sufficient conditions in order to guarantee that a Riemannian manifold is hyperbolic.
See also \cite{T} for the case of Riemannian surfaces.

\section{Main results}

In this section we prove that having a pole is a necessary condition for a hyperbolic uniform graph or a complete Riemannian manifold with bounded local geometry satisfying isoperimetric inequality. Therefore, it can be removed as hypothesis in the statements of Theorems \ref{t:iigraph} and \ref{th: isoperimetric_manifolds}.

In fact, we are going to prove the following results.

\begin{theorem} \label{t:iigraph2}
Given a hyperbolic uniform graph $\Gamma$, then $h(\Gamma)>0$ if and only if $\partial_\infty \Gamma$ is uniformly perfect for some visual metric and $\Gamma$ is an infinite graph with a pole.
\end{theorem}

\begin{theorem} \label{th: isoperimetric_manifolds2}
Let $X$ be a hyperbolic complete Riemannian manifold with bounded local geometry.
Then, $h(X)>0$ if and only if $\partial_\infty X$ is uniformly perfect and $X$ is non-compact and has a pole.
\end{theorem}

These results are a consequence of Theorems \ref{t:iigraph}, \ref{th: isoperimetric_manifolds}, and the two following results.

\begin{theorem} \label{t:graphpole}
Let $\Gamma$ be a uniform graph. If $\Gamma$ is hyperbolic and $h(\Gamma)>0$, then $\Gamma$ is an infinite graph with a pole.
\end{theorem}

\begin{theorem} \label{th:manifoldspole}
Let $X$ be a complete Riemannian manifold with bounded local geometry.
If $X$ is hyperbolic and $h(X)>0$, then $X$ is non-compact and has a pole.
\end{theorem}

\medskip

Let us now proceed with the proof of Theorem \ref{t:graphpole}.

\begin{proof}
Since $h(\Gamma)>0$, we have that $\Gamma$ is an infinite graph.

Seeking for a contradiction assume that $\Gamma$ does not have a pole.
Let $\mu$ be a constant such that $\G$ is $\mu$-uniform.

Fix $v\in V(\G)$ and denote by $K$ the union of the geodesic rays starting from $v$.
Since $\Gamma$ does not have a pole, for each $n$ there exists $v_n\in V(\G)$ with $d_\G(v_n,K) \ge 4n$.
Let $\eta_n : [0,\ell_n] \rightarrow \G$ be a geodesic joining $v_n$ with $K$ and such that $\eta_n(0)= v_n$ and
$d_\G(v_n,K) = L(\eta_n) = \ell_n  \ge 4n$.

Given $s \in \RR$, denote by $\lfloor s \rfloor$ the lower integer part of $s$, i.e., the largest integer not greater than $s$.
Assume that the ball $B_\G \big(\eta_n(t),\lfloor \d_{th}(\G) \rfloor+1 \big)$ intersects every geodesic joining $v_n$ with some point of $K$,
for some $t \in \ZZ^+$ with $\lfloor \d_{th}(\G) \rfloor  < t \le 4n \le \ell_n$.
Thus, the connected component $A_n$ of $\G \setminus B_\G \big(\eta_n(t),\lfloor \d_{th}(\G) \rfloor+1 \big)$ containing $v_n$ satisfies
that its boundary $\p A_n$ is contained in the sphere $S_\G \big(\eta_n(t),\lfloor \d_{th}(\G) \rfloor \big)$.

Since $\G$ is $\mu$-uniform,
$$
| \p A_n |
\le \big| S_\G \big(\eta_n(t),\lfloor \d_{th}(\G) \rfloor \big) \big|
\le \mu^{\lfloor \d_{th}(\G) \rfloor} .
$$
Thus
$$
t - \lfloor \d_{th}(\G) \rfloor
= \big| V(\G) \cap \eta_n \big(\big[ 0, t - \lfloor \d_{th}(\G) \rfloor - 1 \big] \big) \big|
\le | A_n |
\le h(\G)^{-1} | \p A_n |
\le h(\G)^{-1} \mu^{\lfloor \d_{th}(\G) \rfloor} ,
$$
and we conclude that
$$
t \le \lfloor \d_{th}(\G) \rfloor + h(\G)^{-1} \mu^{\lfloor \d_{th}(\G) \rfloor}
=: M .
$$

Hence, if $2n> M$, there exists a geodesic $\s_n$ from $v_n$ to $K$ such that
$\s_n \cap B_\G \big(\eta_n(2n),\lfloor \d_{th}(\G) \rfloor+1 \big) = \emptyset$.
Note that $d_\G(\eta_n(2n),K) = \ell_n - 2n \ge 2n$.
Let $x_n$ and $y_n$ be the endpoints of $\eta_n$ and $\s_n$ in $K$, respectively,
and consider the geodesic triangle $T^n=\{v_n,x_n,y_n\}$ in $\G$.
Since $\d_{th}(\G) < \lfloor \d_{th}(\G) \rfloor+1$,
we have
$$
d_\G(\eta_n(2n),\s_n)
\ge \lfloor \d_{th}(\G) \rfloor+1
> \d_{th}(\G),
$$
and so, there exists $z_n \in [x_ny_n]$ with $d_\G(\eta_n(2n),z_n) \le \d_{th}(\G)$ (see Figure \ref{Fig:pole}).

\begin{figure}[h]
\centering
\includegraphics[scale=0.4]{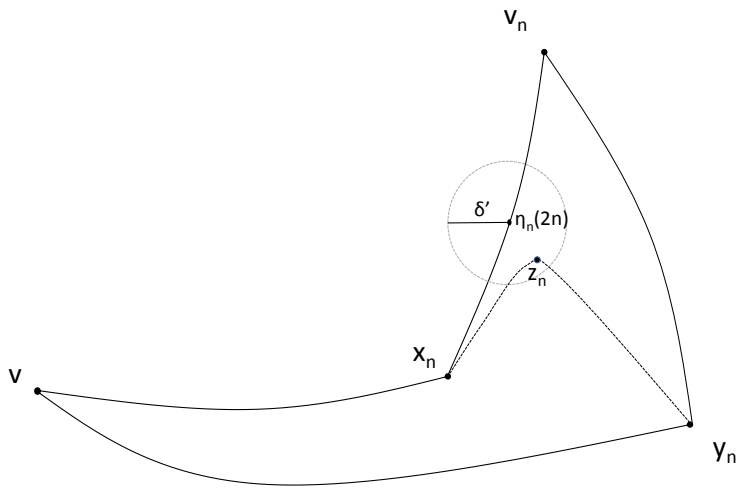}
\caption{If there is a geodesic $\s_n=[v_ny_n]$ from $v_n$ to $K$ which does not intersect the ball  $B_\G \big(\eta_n(2n),\d' \big)$, where $\d':=\lfloor \d_{th}(\G) \rfloor+1$, then there is a point $z_n\in [x_ny_n]$ which is far from~$K$.}
\label{Fig:pole}
\end{figure}

Consider now the geodesic triangle $T_n=\{v,x_n,y_n\}$ in $\G$, and $z_n \in [x_ny_n]$.
Since $[vx_n]\cup [vy_n] \subset K$,
we have
$$
\d_{th}(\G)
\ge d_\G(z_n,[vx_n]\cup [vy_n])
\ge d_\G(z_n,K)
\ge d_\G(\eta_n(2n),K) - d_\G(\eta_n(2n),z_n)
\ge 2n - \d_{th}(\G),
$$
and so, $\d_{th}(\G) \ge n$ for every $n$.
This is the contradiction we were looking for, since $\G$ is hyperbolic, and we conclude that $\Gamma$ has a pole.
\end{proof}

\medskip

Finally, let us now proceed with the proof of Theorem \ref{th:manifoldspole}.

\begin{proof}
Since $h(X)>0$, we have that $X$ is non-compact.

Let $\Gamma$ be an $\varepsilon$-net in $X$.
Since $h(X)>0$, Proposition \ref{Prop: Kanai_net} gives that $h(\Gamma)>0$.

Note that $X$ is a proper geodesic space since it is a complete Riemannian manifold.
Also, $\G$ is a proper geodesic space since it is a uniform graph by Lemma \ref{lema: uniform-net}.

By Proposition \ref{Prop: net_qi}, $X$ and $\Gamma$ are quasi-isometric and so, $\Gamma$ is hyperbolic by Theorem \ref{th: stability_hyp}.
Theorem \ref{t:graphpole} gives that $\Gamma$ has a pole, and
by Proposition \ref{Prop: qi-pole}, $X$ has a pole.
\end{proof}

\begin{remark} Theorem \ref{t:iigraph2} also allows to improve Theorem 4.12 in \cite{MR2} which is based in the same results. In the same way, having a pole is a necessary condition and can be removed as hypothesis.
\end{remark}

\begin{remark} Theorem 5.11 in \cite{Me} states that for geodesic, visual, Gromov hyperbolic spaces, having uniformly perfect boundary at infinity is equivalent to being \emph{uniformly equilateral}. Thus, theorems \ref{t:iigraph2} and \ref{th: isoperimetric_manifolds2} can also be written using this property and without using the boundary at infinity.
\end{remark}

\

\noindent{\bf Acknowledgments:}
We would like to thank to Paloma Mart\'{\i}n for her support during this research.
This work was supported in part by two grants
from Ministerio de Econom{\'\i}a y Competitividad, Agencia Estatal de
Investigaci\'on (AEI) and Fondo Europeo de Desarrollo Regional (FEDER) (MTM2015-63612P, MTM2016-78227-C2-1-P and MTM2017-90584-REDT), Spain.


\begin{thebibliography}{99}


\bibitem{AAD} Abu-Ata, M., Dragan, F. F.: Metric tree-like structures in real-life networks: an empirical study. {\it Networks} {\bf 67}, 49-68 (2016)

\bibitem{ASM} Adcock, A. B., Sullivan, B. D. and Mahoney, M. W.: Tree-like
structure in large social and information networks. In: 13th Int
Conference Data Mining (ICDM), pp. 1-10. IEEE, Dallas, Texas, USA (2013)

\bibitem{ABCD} Alonso, J., Brady, T., Cooper, D., Delzant, T., Ferlini, V., Lustig, M., Mihalik, M., Shapiro, M. and Short, H.: Notes on word hyperbolic groups. In: E. Ghys, A. Haefliger, A. Verjovsky (Eds.),
Group Theory from a Geometrical Viewpoint, World Scientific, Singapore (1992)

\bibitem{APR} Alvarez, V., Pestana, D., Rodr\'{\i}guez, J. M.: Isoperimetric inequalities in Riemann surfaces of infinite type.
{\it Rev. Mat. Iberoamericana} {\bf 15}, 353-427 (1999)

\bibitem{A1} Ancona, A.: Negatively curved manifolds, elliptic operators, and Martin boundary.
{\it Annals of Math.} {\bf 125}, 495-536 (1987)

\bibitem{A2} Ancona, A.: Positive harmonic functions and hyperbolicity. In: Potential Theory, Surveys
and Problems, eds. J. Kr\'al et al., Lecture Notes in Math., No. 1344, pp. 1-24. Springer-Verlag, Heidelberg (1988)

\bibitem{A3} Ancona, A.: Theorie du potentiel sur les graphes et les varieties. In: Ecol\'e d'Et\'e de
Probabilit\'es de Saint-Flour XVII-1988, eds. A.Ancona et al., Lecture Notes in Math.,
No. 1427, Springer-Verlag, Heidelberg (1990)

\bibitem{AC} Anderson, M., Cheeger, J.: $C^\a$-compactness for manifolds with Ricci curvature and
injectivity radius bounded below. {\it J. Diff. Geom.} {\bf 35}, 265-281 (1992)


\bibitem{BRS} Bermudo, S., Rodr\'{\i}guez, J. M. and Sigarreta, J. M.:
Computing the hyperbolicity constant. {\it Comput. Math. Appl.} {\bf 62}, 4592-4595 (2011)


\bibitem{BJ} Bishop, C. J., Jones, P. W.: Hausdorff dimension and Kleinian groups.
{\it Acta Math.} {\bf 179}, 1-39 (1997)

\bibitem{BH} Bridson, M., Haefliger, A.: Metric spaces of non-positive curvature.
Springer-Verlag, Berlin (1999)

\bibitem{BHB1} Brinkmann, G., Koolen J. and Moulton ,V.: On the hyperbolicity of chordal
graphs. {\it Ann. Comb.} {\bf 5}, 61-69 (2001)



\bibitem{BBI} Burago, D., Burago, Y., Ivanov, S.: A course in metric geometry. Graduate
Studies in Mathematics, 33, AMS, Providence, RI (2001)

\bibitem{Bu1} Buser, P.: A note on the isoperimetric constant.
{\it Ann. Sci. \'Ecole Normale Sup.} {\bf 15}, 213-230 (1982)

\bibitem{Bu} Buser, P.: Geometry and Spectra of Compact Riemann
Surfaces. Birkh\"auser, Boston (1992)

\bibitem{BS}  Buyalo, S.,  Schroeder, V.: \emph{Elements of Asymptotic
Geometry.} EMS Monographs in Mathematics, Germany (2007)

\bibitem{CGPR} Cant\'on, A., Granados, A., Portilla, A., Rodr\'{\i}guez, J. M.:
Quasi-isometries and isoperimetric inequalities in planar domains.
{\it J. Math. Soc. Japan} {\bf 67}, 127-157 (2015)

\bibitem{Cao} Cao, J.: Cheeger isoperimetric constants of Gromov-hyperbolic spaces with quasi-pole.
{\it Commun. Contemp. Math.} {\bf 4}(2), 511-533 (2000)


\bibitem{C2} Chavel, I.:
Isoperimetric inequalities: differential geometric and analytic perspectives.
Cambridge University Press, Cambridge (2001)

\bibitem{Ch} Cheeger, J.: A lower bound for the smallest eigenvalue of
the Laplacian. In: {\it Problems in Analysis}, pp. 195-199. Princeton University
Press, Princeton (1970)

\bibitem{CYY} Chen, B., Yau, S.-T. and Yeh, Y.-N.: Graph homotopy and Graham homotopy.
{\it Discrete Math.} {\bf 241}, 153-170 (2001)


\bibitem{CMN} Clauset, A., Moore, C., Newman, M. E. J.: Hierarchical structure and
the prediction of missing links in networks. {\it Nature} {\bf 453}, 98-101 (2008)


\bibitem{DHHKMW} Dress, A., Holland, B., Huber, K.T., Koolen, J.H., Moulton, V., Weyer-Menkhoff, J.:
$\Delta$ additive and $\Delta$ ultra-additive maps, Gromov's trees, and the Farris transform. {\it Discrete Appl. Math.} {\bf 146}(1), 51-73 (2005)

\bibitem{DMT} Dress, A., Moulton, V., Terhalle, W.: T-theory: an overview. {\it Europ. J. Combin.} {\bf 17}, 161-175 (1996)

\bibitem{FM1} Fern\'andez, J. L., Meli\'an, M. V.:
Bounded geodesics of Riemann surfaces and hyperbolic manifolds.
{\it Trans. Amer. Math. Soc.} {\bf 347}, 3533-3549 (1995)

\bibitem{FM2} Fern\'andez, J. L., Meli\'an, M. V.:
Escaping geodesics of Riemannian surfaces. {\it Acta Math.} {\bf 187}, 213-236 (2001)

\bibitem{FMP1} Fern\'andez, J. L., Meli\'an, M. V., Pestana, D.:
Quantitative mixing results and inner functions. {\it Math. Ann.} {\bf 337}, 233-251 (2007)

\bibitem{FMP2} Fern\'andez, J. L., Meli\'an, M. V., Pestana, D.:
Expanding maps, shrinking targets and hitting times.
{\it Nonlinearity} {\bf 25}, 2443-2471 (2012)

\bibitem{FR1} Fern\'andez, J. L., Rodr{\'\i}guez, J. M.:
The exponent of convergence of Riemann surfaces: Bass Riemann surfaces.
{\it Annales Acad. Sci. Fenn. A. I.} {\bf 15}, 165-183 (1990)

\bibitem{GH} Ghys, E., de la Harpe, P.:
{\it Sur les Groupes Hyperboliques d'apr\`es Mikhael Gromov.}
Progress in Mathematics, Volume 83. Birkh\"auser, Berlin (1990)

\bibitem{GPPRT} Granados, A., Pestana, D., Portilla, A., Rodr\'{\i}guez, J. M., Tour\'{\i}s, E.:
Stability of the injectivity radius under quasi-isometries and applications to isoperimetric inequalities.
{\it RACSAM} {\bf 112}, 1225-1247 (2018)

\bibitem{G1} Gromov, M.: Hyperbolic groups. In: Essays in group theory.
Edited by S. M. Gersten, M. S. R. I. Publ. {\bf 8}, pp. 75-263. Springer, Heidelberg (1987)



\bibitem{K21} Jonckheere, E. A.: Contr\^ole du traffic sur les r\'eseaux \`a
g\'eom\'etrie hyperbolique--Vers une th\'eorie g\'eom\'etrique de la s\'ecurit\'e
l'acheminement de l'information. {\it J. Europ. Syst. Autom.} {\bf 8}, 45-60 (2002)

\bibitem{K22} Jonckheere, E. A. and Lohsoonthorn, P.: Geometry of network security.
{\it Amer. Control Conf.} {\bf ACC}, 111-151 (2004)

\bibitem{K} Kanai, M.: Rough isometries and combinatorial approximations of geometries of noncompact
Riemannian manifolds. {\it J. Math. Soc. Japan} {\bf 37}, 391-413 (1985)



\bibitem{KPKVB} Krioukov, D., Papadopoulos, F., Kitsak, M., Vahdat, A., Bogu\~na, M.:
Hyperbolic geometry of complex networks. {\it Phys. Rev. E} {\bf 82}(3), $\#$036106 (2010)


\bibitem{M} Mart\'{\i}nez-P\'erez, A.:
Chordality properties and hyperbolicity on graphs. \emph{Electron. J. Comb.} \textbf{23}(3), \# P3.51 (2016)

\bibitem{M2} Mart\'{\i}nez-P\'erez, A.: Generalized Chordality, Vertex Separators and Hyperbolicity on Graphs. \emph{Symmetry}, \textbf{9}(10), 199 (2017)


\bibitem{MR1} Mart\'inez-P\'erez, A. Rodr\'iguez J. M.: Cheeger isoperimetric constant of Gromov hyperbolic manifolds and graphs. {\it Commun. Contemp. Math.} \textbf{20}(5) (2018). https://doi.org/10.1142/S021919971750050X

\bibitem{MR2} Mart\'inez-P\'erez, A. Rodr\'iguez J. M.: Isoperimetric inequalities in Riemann surfaces and graphs. Submitted.

\bibitem{MRT} Meli\'an, M. V., Rodr{\'\i}guez, J. M., Tour{\'\i}s, E.:
Escaping geodesics in Riemannian surfaces with pinched variable negative curvature.
Submitted.

\bibitem{Me} Meyer, J.: Uniformly perfect boundaries of Gromov hyperbolic spaces. Ph. D. thesis, University of Zurich (2009)

\bibitem{MoSoVi} Montgolfier, F., Soto, M. and Viennot, L., Treewidth and Hyperbolicity of the Internet. In: 10th IEEE
International Symposium on Network Computing and Applications (NCA), pp. 25–32 (2011)

\bibitem{P} Paulin, F.: On the critical exponent of a discrete group of hyperbolic isometries.
\emph{Differ. Geom. Appl.} {\bf 7}, 231-236 (1997)

\bibitem{PRSV} Pestana, D., Rodr\'{\i}guez, J. M., Sigarreta, J. M., Villeta, M.: Gromov hyperbolic cubic graphs.
{\it Cent. Eur. J. Math.} {\bf 10}, 1141-1151 (2012)

\bibitem{Po} P\'olya, G.: Isoperimetric inequalities in mathematical physics. Princeton University Press (1951)




\bibitem{RT2} Rodr\'{\i}guez, J. M., Tour\'{\i}s, E.: A new characterization of Gromov hyperbolicity for negatively curved surfaces.
{\it Publ. Mat.} {\bf 50}, 249-278 (2006)

\bibitem{Sha1} Shang, Y.: Lack of Gromov-hyperbolicity in colored random networks. {\it Pan-American Math. J.} {\bf 21}(1), 27-36 (2011)

\bibitem{Sha2} Shang, Y.: Lack of Gromov-hyperbolicity in small-world networks. {\it Cent. Eur. J. Math.} {\bf 10}(3), 1152-1158 (2012)

\bibitem{Sha3} Shang, Y.: Non-hyperbolicity of random graphs with given expected degrees. {\it Stoch. Models} {\bf 29}(4), 451-462 (2013)

\bibitem{S} Sullivan, D.: Related aspects of positivity in Riemannian
geometry. {\it J. Diff. Geom.} {\bf 25}, 327-351 (1987)

\bibitem{T} Tour{\'\i}s, E.: Graphs and Gromov hyperbolicity of non-constant negatively curved surfaces.
{\it J. Math. Anal. Appl.} {\bf 380}, 865-881 (2011)

\bibitem{V} V\"ais\"al\"a, J.: Hyperbolic and uniform domains in Banach spaces.
{\it Exposit. Math.} {\bf 23}(3), 187-231 (2005)

\bibitem{WZ} Wu, Y. and Zhang, C.: Chordality and hyperbolicity of a graph. {\it Electr. J. Comb.} {\bf 18}, P43 (2011)

\end{thebibliography}
\end{document}